\documentstyle[12pt]{article}

\def\R{{\hat{R}_{12}}}
\def\ba{\begin{eqnarray}}
\def\ea{\end{eqnarray}}
\def\lb{\label}
\def\be{\begin{equation}}
\def\ee{\end{equation}}

\begin{document}

\title{$SL_{q}(N)$ DIFFERENTIAL CALCULUS FROM THE DIFFERENTIAL CALCULUS
ON $GL_{q}(N)$\thanks{In Proceedings of the VII International Conference
'Symmetry Methods in Physics' (Dubna, July 10-16, 1995)}}

\author{  A.P. ISAEV\thanks{e-mail address: isaevap@thsun1.jinr.dubna.su}
\\
\it Bogoliubov Laboratory of Theoretical Physics \\
\it Joint Institute for Nuclear Research \\
\it 141980 Dubna, Moscow Region, RUSSIA.
}
\date{}
\maketitle

\begin{abstract}
We show that the Faddeev-Pyatov $SL_{q}(N)$ differential algebra \cite{FP}
is a subalgebra of the $GL_{q}(N)$ differential algebra constructed by
Schupp, Watts and Zumino \cite{SWZ}.
\end{abstract}

In the paper \cite{FP}, the bicovariant differential algebra
on $SL_{q}(N)$ (see below (\ref{s2'}), (\ref{s2''}))
with generators $\{ T_{ij}, \; L_{ij}, \;
\tilde{\Omega}_{ij} \}$ $i,j = 1, \dots N$
has been constructed. The
elements $T_{ij}$ are generators of
the algebra $Fun(SL_{q}(N))$, while $L_{ij}$  generate a matrix
of the right-invariant Lie derivatives on $SL_{q}(N)$ and the elements
$\tilde{\Omega}_{ij}$ define the basis of the right-invariant differential
1-forms on $SL_{q}(N)$. It has been shown \cite{FP} that
this algebra is consistent with imposing the conditions:
\be
\lb{isa1}
det_{q}(T) = Det_{q}(L) = 1 \; , \;\; Tr_{q}(\tilde{\Omega}) = 0 \; ,
\ee
where we use
$$
Det_{q}(L) = det_{q}(LT) \frac{1}{det_{q}(T)} \; , \;\;
Tr_{q}(\tilde{\Omega}) = \sum_{i=1}^{N} q^{-N -1 +2i}
\tilde{\Omega}_{ii} \; .
$$
The last condition in (\ref{isa1}) shows that the number of independent
differential 1-forms on $SL_{q}(N)$, in the Faddeev-Pyatov approach,
is equal to $(N^{2}-1)$ just as in the
classical case.

Another bicovariant differential algebra (which was considered as a
differential algebra on $GL_{q}(N)$)
has been proposed in \cite{SWZ} in the framework of the
$R$- matrix approach \cite{FRT} to the
general theory of differential calculi on quantum groups \cite{Wor}.
The defining relations for this algebra have the form:
\be
\lb{zum}
\left\{
\begin{array}{l}
\R \, T_{1} \, T_{2}  =  T_{1} \, T_{2} \, \R \; , \;\;
T_{1} \,  L_{2} = \R \, L_{1} \, \R \, T_{1} \; , \\ \\
T_{1} \, \Omega_{2} =  \R^{-1} \, \Omega_{1} \, \R^{-1} \, T_{1}  \; , \;\;
\R^{-1} \, \Omega_{1} \, \R^{-1} \Omega_{1}  +
\Omega_{1} \, \R^{-1} \, \Omega_{1} \, \R = 0 \; , \\ \\
\R \, L_{1} \, \R \,  L_{1} =
 L_{1} \, \R \,  L_{1} \, \R  \; , \; \;
\R \,  L_{1} \, \R \, \Omega_{1} =
\Omega_{1} \, \R \, L_{1} \, \R  \; ,
\end{array}
\right.
\ee

\be
\lb{8.61}
\left\{
\begin{array}{l}
T_{1} \, \Im_{2} = \R \, \Im_{1} \, \R \, T_{1}  \; , \;\;
 \Im_{1} \, \R^{-1} \, \Omega_{1} \, \R^{-1}  +
\R^{-1} \, \Omega_{1} \, \R^{-1} \, \Im_{1} =  - \R^{-1} \; , \\ \\
\R \,  L_{1} \, \R \, \Im_{1} =
\Im_{1} \, \R \,  L_{1} \, \R   \; , \;\;
\R \, \Im_{1} \, \R \, \Im_{1}  +
\Im_{1} \, \R \, \Im_{1} \, \R^{-1} = 0 \; .
\end{array}
\right.
\ee
Here we rewrite the
defining relations for the differential algebra of \cite{SWZ} in terms
of the right-invariant generators,
use the notation $L$ instead of $Y$, change
$GL_{q}(N)$ $R$-matrix $\R$ to $\R^{-1}$ and
denote right-invariant inner derivations by $\Im_{ij}$.
For the algebra (\ref{zum}), (\ref{8.61}) one can construct
a 1- parametric family of differentials $d_{x}$ satisfying
the usual Leibnitz rule:
\be
\lb{isa2}
d_{x} A = [ \xi_{x}, \, A ]_{\pm} \; , \;\;
\xi_{x} = \frac{q^{N}}{\lambda} \,
Tr_{q} \left( \Omega (1 + x \, W \, L) \right) \; , \;\;
\xi_{x}^{2} =0 \; ,
\ee
where $\lambda = q-q^{-1}$, $W = L(1 - \lambda \Omega \Im)L^{-1}$
(the operator $(1-WL)/\lambda$ defines right-invariant vector
fields on $GL_{q}(N)$), $A$ is a function
of $\{ T, \, \Omega \}$ and $x$ is a parameter.
The last identity in (\ref{isa2}) follows from the relation
(for all $x$):
$$
\R \Omega_{x} \R^{-1} \Omega_{x} =
- \Omega_{x} \R^{-1} \Omega_{x} \R^{-1} \; , \;\;
\Omega_{x} = \Omega \, (1 + x W L) \; .
$$
In the paper \cite{SWZ}, there was suggested a way to
reduce the $GL_{q}(N)$ algebra (\ref{zum}), (\ref{8.61})
to an algebra which was interpreted as
a differential algebra on $SL_{q}(N)$. Unfortunately, this
reduction to the $SL_{q}(N)$ case leads to
$N^{2}$ independent differential 1-forms $\Omega_{ij}$.
Below we demonstrate
that the $SL_{q}(N)$ differential algebra \cite{FP},
with a correct number $(N^{2}-1)$
of independent differential 1-forms, can be obtained as a subalgebra
of the $GL_{q}(N)$ differential algebra (\ref{zum}), (\ref{8.61}).

Let us pass in eqs. (\ref{zum}), (\ref{8.61})
to a new basis of differential forms, inner derivations
and new generators of the quantum group $Fun(GL_{q}(N))$ (see \cite{Is}):
\be
\lb{s1}
\Omega \rightarrow \Omega^{L} = L \, \Omega  \; , \;\;
\Im \rightarrow \Im^{L} =  \Im \, L^{-1} \; , \;\;
T \rightarrow (det_{q} T)^{-1/N} \, T \; .
\ee
Then, relations (\ref{zum}), (\ref{8.61})
take the following form in terms of the new generators:
\be
\lb{s2}
\left\{
\begin{array}{c}
\R \, T_{1} \, T_{2}  =  T_{1} \, T_{2} \, \R \; , \;\;
 q^{2/N}  T_{1} \, L_{2} =  \R \,  L_{1} \, \R \, T_{1} \; , \;\; \\ \\
T_{1} \, \Omega^{L}_{2} =  \R \, \Omega^{L}_{1} \, \R^{-1} \, T_{1}  \; , \;\;
\R \, \Omega^{L}_{1} \, \R \, \Omega^{L}_{1}  +
\Omega^{L}_{1} \, \R \, \Omega^{L}_{1} \, \R^{-1} = 0 \; ,  \\ \\
\R \,  L_{1} \, \R \,  L_{1} =
L_{1} \, \R \, L_{1} \, \R   \; , \;\;
\R \, L_{1} \, \R \, \Omega^{L}_{1} =
\Omega^{L}_{1} \, \R \, L_{1} \, \R  \; ,
\end{array}
\right.
\ee

\be
\lb{s3}
\left\{
\begin{array}{l}
T_{1} \, \Im^{L}_{2} =  \R \, \Im^{L}_{1} \, \R^{-1} \, T_{1} \; , \;\;
\R \, {\Im^{L}}_{1} \, \R^{-1} {\Im^{L}}_{1}  +
{\Im^{L}}_{1} \, \R^{-1} \, {\Im^{L}}_{1} \, \R^{-1} = 0 \; , \\ \\
\R \,  L_{1} \, \R \, {\Im^{L}}_{1} =
{\Im^{L}}_{1} \, \R \,  L_{1} \, \R   \; , \;\;
\end{array}
\right.
\ee

\be
\lb{ss3}
{\Im^{L}}_{1} \, \R^{-1} \, \Omega^{L}_{1} \, \R +
\R \, \Omega^{L}_{1} \, \R^{-1} \, {\Im^{L}}_{1} = -  \R \; .
\ee
Now we show that the subalgebra
(\ref{s2}) contains the Faddeev-Pyatov
$SL_{q}(N)$ subalgebra $\{ T, \; L, \; \tilde{\Omega} \}$
where we have introduced the traceless
generators $\tilde{\Omega}$ defined by
\be
\lb{ss}
\tilde \Omega = \Omega^{L} -
{1 \over N_q } Tr_q \Omega^{L} \cdot {\bf 1} \; ,
\; \; N_{q} = \frac{q^{N} - q^{-N}}{q-q^{-1}} \; .
\ee
Indeed, after substituting (\ref{ss}) into relations (\ref{s2})
and using formulas (which are a consequence of the fourth equation
in (\ref{s2}))
$$
(Tr_{q}\Omega^{L})^{2} = 0 \; , \;\;
[Tr_{q} \Omega^{L} , \; \tilde{\Omega} ]_{+} =
\lambda q^{N} (\kappa_{q} -1) \tilde{\Omega}^{2} \; ,
$$
where $\kappa_{q} = \lambda q^{N} (N_{q} + \lambda q^{N})^{-1}$,
we deduce the algebra \cite{FP}:
\be
\lb{s2'}
\left\{
\begin{array}{c}
\R \, T_{1} \, T_{2}  =  T_{1} \, T_{2} \, \R \; , \;\;
 q^{2/N}  T_{1} \, L_{2} =  \R \,  L_{1} \, \R \, T_{1} \; , \;\; \\ \\
T_{1} \, \tilde{\Omega}_{2} =  \R \,
\tilde{\Omega}_{1} \, \R^{-1} \, T_{1}  \; , \;\;      \\ \\
\R \,  L_{1} \, \R \,  L_{1} =
L_{1} \, \R \, L_{1} \, \R   \; , \;\;
\R \, L_{1} \, \R \, \tilde{\Omega}_{1} =
\tilde{\Omega}_{1} \, \R \, L_{1} \, \R  \; ,
\end{array}
\right.
\ee

\be
\lb{s2''}
\R \, \tilde{\Omega}_{1} \, \R \, \tilde{\Omega}_{1}  +
\tilde{\Omega}_{1} \, \R \, \tilde{\Omega}_{1} \, \R^{-1} =
\kappa_{q} (\tilde{\Omega}^{2} + \R \tilde{\Omega}^{2} \R) \; .
\ee
The algebra of differential 1-forms, with defining relations
(\ref{s2''}), has been proposed in \cite{IP}.
It is clear that the whole algebra (\ref{s2})-(\ref{ss3}) includes
the subalgebra $\{ T, \; L, \; \tilde{\Omega}, \; \Im \}$ with
defining relations
(\ref{s2'}), (\ref{s2''}), (\ref{s3}) and,
instead of the relation (\ref{ss3}), we also have to take
\be
\lb{ss4}
{\Im^{L}}_{1} \, \R^{-1} \, \tilde{\Omega}_{1} \, \R +
\R \, \tilde{\Omega}_{1} \, \R^{-1} \, {\Im^{L}}_{1} =
\frac{\kappa_{q}}{\lambda(1-\kappa_{q})} -  \R \; .
\ee

One can show that the algebra (\ref{s2})-(\ref{ss3})
(and, therefore, the subalgebra
(\ref{s2'}), (\ref{s2''}), (\ref{s3}), (\ref{ss4}))
possesses an additional central element \cite{Is}:
\be
\lb{s5}
Z = Det_{q} \left( (\overline{W} \, W)^{-1} \right) \equiv
det_{q} \left( (\overline{W} \, W)^{-1} \, T \right)
\, \frac{1}{det_{q}(T)}  \; , \;\;
\ee
$$
Z = Det_{q}(W^{-1}) \cdot Det_{q}(\overline{W}^{-1}) \; ,
$$
where the operators $\overline{W} = 1 - \lambda \Im^{L} \Omega^{L}$ and
$W = 1 - \lambda \Omega^{L} \Im^{L}$ satisfy the commutation relations:
\be
\lb{s7}
\begin{array}{c}
\R^{-1} \, \overline{W}_{1} \, \R \, W_{1} =
W_{1} \, \R^{-1} \,  \overline{W}_{1} \, \R \; , \;\;      \\ \\
\R^{-1} \, \overline{W}_{1} \, \R^{-1} \, \overline{W}_{1} =
\overline{W}_{1} \, \R^{-1} \,  \overline{W}_{1} \, \R^{-1} \; , \;\; \\ \\
\R^{-1} \, {W}_{1} \, \R^{-1} \, {W}_{1} =
{W}_{1} \, \R^{-1} \,  {W}_{1} \, \R^{-1} \; . \;\;
\end{array}
\ee
Note that from these relations we immediately (see \cite{Maj}) obtain:
\be
\lb{s8}
\R^{-1} \, (\overline{W}\, W)_{1} \, \R^{-1} \, (\overline{W}\, W)_{1} =
(\overline{W}\, W)_{1} \, \R^{-1} \,  (\overline{W}\, W)_{1} \, \R^{-1} \; .
\ee
We stress also that the combination
$(\overline{W}\, W)$ is independent of the
scalar generator $Tr_{q} \Omega^{L}$:
$$
(\overline{W}\, W) = \frac{1}{1 - \kappa_{q}}
- \lambda [ \tilde{\Omega} , \, \Im^{L} ]_{+} +
\lambda^{2} \, (1-\kappa_{q})
\, \Im^{L} \, \tilde{\Omega}^{2} \, \Im^{L} \; .
$$
It means that if we fix the central element (\ref{s5}), e.g.
in the form $Z=1$, we, therefore, impose a new additional relation
on the generators $\Im^{L}_{ij}$. One can try to consider the
relation $Z=1$ as a constraint which
finally reduces the
subalgebra generated by $\{ T, \; L, \; \tilde{\Omega}, \; \Im \}$
to the whole differential algebra on $SL_{q}(N)$. However, it is
not quite clear: is the relation $Z=1$ sufficient for removing
the redundant generator $Tr_{q} \Im^{L}$ from
this algebra to obtain the complete $SL_{q}(N)$ reduction?

Unfortunately, neither of the differential operators (\ref{isa2})
(for the Schupp-Watts-Zumino algebra (\ref{zum}), (\ref{8.61}))
yields an apropriate differential $d$ (satisfying the usual Leibnitz rule)
for the Faddeev-Pyatov differential
subalgebra (\ref{s2'}), (\ref{s2''}).
Thus, up to now, we have only one self-consistent construction
of the exterior differential $d$ (for the algebra
(\ref{s2'}), (\ref{s2''})) which satisfies the deformed
Leibnitz rule (see \cite{FP}). Nevertheless, we hope
that the method of obtaining the
differential algebra on $SL_{q}(N)$
(\ref{s2'}), (\ref{s2''}) as a subalgebra
of the differential
algebra on $GL_{q}(N)$ (\ref{zum}), (\ref{8.61}) will help us
to clarify the Hopf and star structure of this
subalgebra and also may help us to construct
the explicit realization
(e.g. via deformed commutator with nilpotent
generators (\ref{isa2})) of an exterior differential $d$
obeying the deformed Leibnitz rule \cite{FP}.

To conclude this report I would like to stress that the situation with
bicovariant differential calculi on $SO_{q}(N)$ and $Sp_{q}(2n)$
groups is obscured up to now. Woronowicz's approach \cite{Wor} leads,
in this case, to the differential algebras which are not of the
Poincar\'{e}- Birkhoff- Witt type and all attempts to improve
this situation were failed (for the discussion see \cite{AIP}).

{\bf Acknowledegment}  \par
I would like to thank G.E.Arutyunov, P.N.Pyatov and
A.A.Vladimirov for many valuable discussions.

This research was supported in part by the
Russian Foundation of Fundamental Research
(grant N\underline{o} 95-02-05679-a), by the
International Science Foundation (grant RFF 300)
and by the INTAS grant
N\underline{o} 93-127.

\end{document}